\documentclass[runningheads]{llncs}
\usepackage{graphicx}
\usepackage{url}
\begin{document}
\title{A Composite Index Method for Optimization Benchmarking with Application to the Backhaul Profit Maximization Problem}
%
%
\author{Yulan Bai\orcidID{0000-0003-0305-3274} \and
Eli Olinick\orcidID{0000-0001-7856-1826} }
\authorrunning{Y. Bai and E. Olinick}
%
\institute{Department of Operations Research and Engineering Management, Southern Methodist University, Dallas TX 75205, USA 
\\
\email{\{yulanb,olinick\}@smu.edu}}
\maketitle              
\begin{abstract}
We propose a multi-criteria Composite Index Method (CIM) to compare the performance of alternative approaches to solving an optimization problem. The CIM is  convenient in those situations when neither approach dominates the other when tested on different sizes of problem instances. The CIM takes problem  instance size and multiple performance criteria into consideration within a weighting scheme to produce a single number that measures the relative improvement of one alternative over the other. Different weights are given to each dimension based on their relative importance as determined by the end user. We summarize the successful application of the CIM to an ${\cal NP}$-hard combinatorial optimization problem known as the backhaul profit maximization problem (BPMP).  Using the CIM we tested a series of eleven techniques for improving solution time using CPLEX to solve two different BPMP models proposed in the literature.


\end{abstract}
\section{Introduction}

Using solution time as the key performance measure is a long-standing standard practice in the optimization literature. However, now that computing environments take advantage of multiple processors and multiple threads while supporting concurrent running of multiple CPU-intensive processes have become commonplace, measuring solution time is no longer straight-forward.  Furthermore, it is often the case that there is a “crossover point” in problem instance size below which one approach is generally “faster” than another, but above which the second approach is faster.  In this situation the second approach would usually be favored because the emphasis in the literature is on solution time as a function of problem instance size.  In this study, however, we consider the practical question of making a recommendation to a user who frequently solves problems that range in size around the crossover point, and propose a multi-criteria framework for comparing competing solution approaches. We propose a Composite Index Method (CIM) that considers several weighted performance measure factors and calculates a single number (a composite index) to measure the relative performance of two competing solution approaches.  

The CIM was developed to evaluate two proposed mixed integer programming (MIP) formulations of the backhaul profit maximization problem (BPMP), the node-arc and triples formulations, each of which can be enhanced with a variety of solution techniques (e.g., branching-rules and cutting planes). The results of our application of the CIM to the BPMP are discussed in \cite{BPMP}. In this paper we focus on  the process of using the CIM to arrive at the final ``candidate" models in \cite{BPMP}.
This falls into the area of optimization benchmarking.  Beiranvand et al. \cite{beiranvand2017best} provide a recent comprehensive review of the benchmarking literature for optimization problems. As far as we know, the first published study in optimization benchmarking was by Hoffman et al. \cite{hoffman1953computational}, in which different methods were proposed for linear programming and different test instances were used to compare algorithms based on the measures of CPU time, number of iterations, and convergence rate. Another important early paper by Box \cite{box1966comparison} considered the importance of model size and the number of function evaluations during comparison. Later, many researchers explored optimization benchmarking in various applications such as unconstrained optimization, nonlinear least squares, global optimization, and derivative-free optimization. Crowder et al. \cite{crowder1979reporting} proposed standards and guidelines for benchmarking algorithms. According to \cite{beiranvand2017best}, the Performance Profile proposed by Dolan and Mor{\'e} \cite{dolan2002benchmarking} has become the ``gold standard" for optimization benchmarking (over 4,000 citations so far).

Given a set of solution approaches to an optimization problem, the procedure for using the Performance Profile may be summarized as follows. First, a single performance measure is selected (usually the computing time). Second, all candidate solution approaches are applied to each of a set of problem instances and the best-performing approach for each of the tested instances is used as the benchmark for assessing the performance of all of the other candidate approaches on that particular problem instance. In the case of CPU time as the selected metric,  the relative performance of a particular approach on a particular instance is measured by a performance ratio obtained by dividing the CPU time of that particular approach by the CPU time of the best-performing approach. For a particular solution approach, the Performance Profile plots the cumulative distribution function of the performance ratio, the percentage of  instances for which the ratio is less than $x$,  over the range $1 \le x \le \infty$.

In cases such as our BPMP study where the Performance Profiles of the candidate approaches intersect and cross each other, it may be unclear which approach is the best overall. Problem size is an important consideration in the BPMP use case; this makes the Performance Profile inappropriate since it treats each problem instance equally regardless of its size. Furthermore, the Performance Profile only uses a single performance measure, which makes it difficult to use when there are multiple performance criteria such as when comparing the trade-off between solution time and solution quality with heuristics, or executing a solution approach on a system with multiple processors and/or threads.

To the best of our knowledge, parallel computing has received much less attention in the literature on optimization benchmarking despite the fact that it is now widely used in applied optimization. Barr and Hickman did pioneering studies in this area \cite{barr1993reporting,BarrHickman94} and proposed solutions to the challenges parallelization brings to benchmarking. However, they did not suggest using a single measure for easy comparison. Hence, the CIM is an initial step in closing a gap in the literature.

The rest of this paper is structured as follows. We propose our Composite Index Method for benchmarking in Section \ref{sec:CIM}. We describe the BPMP and node-arc and triples formulations in Section \ref{sec:BPMP}. We illustrate the application of the CIM to the node-arc formulation in Section \ref{sec:Node-Arc} and summarize our results from applying the CIM to the triples formulation in Section \ref{sec:Triples}.  We draw conclusions in Section \ref{sec:Conclusion}.

\section{Performance Evaluation Using Composite Index Method (CIM)}
\label{sec:CIM}
We use the CPLEX MIP solver  \cite{CPLEX} to solve the BPMP instances described in \cite{BO,BPMP}.
There are three kinds of “solution time” in the CPLEX output: “CPU time”, “real time”, and “ticks”. CPU time is a measure of the total time used by CPLEX to find an optimal solution; it is the total time used by all threads. Real time (also called wall clock time) is the time that elapsed during the CPLEX run. Both measures can vary noticeably between runs with identical input on identical hardware. Therefore, we solve each problem instance three times in each experiment and report the average CPU and real time over the three runs. The tick metric, also called deterministic time, is a proprietary measure of computation effort based on counting the number of instructions executed by the CPLEX solver and therefore shows no variation between multiple runs with the same inputs on a given hardware configuration \cite{Ticks}. 

For each of the time measures described above, we use a speedup measure to compare the solution time of two solution approaches, approach 1 versus approach 2. Note that in the BPMP application described herein, a solution approach is essentially a MIP model for the BPMP implemented in AMPL \cite{AMPL} and solved with CPLEX. In general, a solution approach could be a combination of a MIP model and an optimization algorithm. Hereinafter, Speedup is defined as the ratio

\[ \mbox{Speedup} = (\mbox{Model 1 solution time}) \div (\mbox{Model 2 solution time}). \]

If Speedup $> 1$, model 2 is solved Speedup times faster than model 1; if Speedup $= 1$ , model 2 has the same solution time as model 1,  and  if Speedup $< 1$, model 1 is solved $1/\mbox{Speedup}$  times faster than model 2.
Due to the fact that CPU and real time are not completely reproducible, we suggest that neither one should be the sole basis for comparing solution approaches. Typically, ticks and real time are positively correlated (as are ticks and CPU time), however there does not appear to be a fixed relationship between ticks and the two time measures. For this reason, we cannot use ticks as the single measure to compare two models either. Instead, we propose a weighted combination of all three time measures. 

For a given problem size, $n$, and timing measure (CPU time, real time, or ticks), we calculate a composite index based on a weighted combination of the minimum, median, mean, and maximum speedups among a set of problem instances. Thus, we obtain three composite indices: $C_n$, $R_n$, $T_n$ for CPU time, real time, and ticks, respectively. To calculate these indices we denote the minimum, mean, median, and maximum speedups in CPU time by $C_{min}$, $C_{mea}$, $C_{med}$, and $C_{max}$, respectively, and define $R_{min}$, $R_{mea}$, $R_{med}$, $R_{max}$, $T_{min}$, $T_{mea}$, $T_{med}$, and $T_{max}$ as the corresponding speedups for real time and ticks. Additionally, we define $\omega_{min}$, $\omega_{mea}$ $\omega_{med}$ and $\omega_{max}$ for  weighting of the minimum, mean, median and maximum statistics. We also define $\bar\omega $ as the summation of $\omega_{min}$, $\omega_{mea}$, $\omega_{med}$ and $\omega_{max}$. The relative weights for CPU, real time, and ticks are $\omega_c$, $\omega_r$, and $\omega_t$, respectively. Using this notation, the three composite indices are calculated as follows:
\begin{eqnarray*}
&& C_n  = (\omega_{min} C_{min} + \omega_{mea} C_{mea} + \omega_{med} C_{med} + \omega_{max} C_{max})/\bar{\omega} \\
&& R_n  = (\omega_{min} R_{min} + \omega_{mea} R_{mea} + \omega_{med} R_{med} + \omega_{max} R_{max})/\bar{\omega}\\
&& T_n  = (\omega_{min} T_{min} + \omega_{mea} T_{mea} + \omega_{med} T_{med} + \omega_{max} T_{max})/\bar{\omega}
\end{eqnarray*}
Next, we calculate a composite index, $I_n$, for problem size $n$ as a weighted combination of indices $C_n$, $R_n$, and $T_n$: 
\[ I_n = (\omega_c C_n + \omega_r R_n + \omega_t T_n)/(\omega_c + \omega_r + \omega_t ).\]

Given a set of problem instance sizes, $\mathcal{S},$ and weight $\omega_s$ for each $s \in \mathcal{S}$, we calculate the grand composite index (GCI), which is the weighted sum of the composite indices for each problem size given by
\[ \mbox{GCI} =( \sum_{s \in \mathcal{S}} \omega_s I_s )/\sum_{s \in \mathcal{S}}\omega_s.  \]

If the grand composite index  $\mbox{GCI} > 1$, we say that model 2 performs GCI times better than model 1; if $\mbox{GCI} = 1$ , model 2 performs the same as model 1,  and  if $\mbox{GCI} < 1$, model 1 performs $1/\mbox{\mbox{GCI}}$  times better than model 2.

In summary, the composite index method (CIM) seeks to find a single  number, GCI, in a parallel computing environment, to decide which solution approach is better, through instance testing. To do so, we first need to decide the performance measures; usually more than one measure is needed. Second, multiple runs are needed to reduce the variance of the measures for the same instance by averaging the measure. Third, for a fixed problem size, multiple instances should be randomly sampled. Along with the mean measure over the different instances of the same problem size, we consider the minimum (min), median, and maximum (max) measure to diminish the effects of outliers. The consideration of min, mean, median and max, is inspired by the famous PERT concept of project management, in which pessimistic, optimistic, and most likely task-completion times are considered with different weights.  Finally, comparisons are made over a range of problem sizes and weighted accordingly.

The steps described above are for comparing two solution approaches. In Section \ref{sec:Node-Arc} we describe how we apply CIM iteratively to compare multiple solution approaches. We illustrate this iterative process by applying it to the BPMP in Sections \ref{sec:Node-Arc} and \ref{sec:Triples}. 

\section{The Backhaul Profit Maximization Problem (BPMP)}
\label{sec:BPMP}
\begin{figure}
\includegraphics[width=\textwidth]{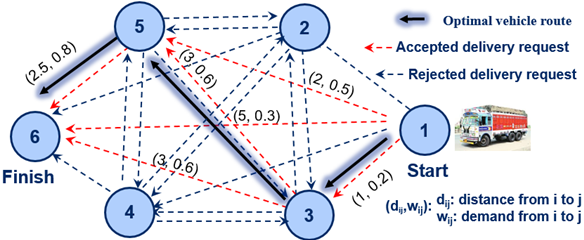}
\caption{BPMP Example.} 
\label{fig:BPMP}
\end{figure}

The BPMP requires simultaneously solving two problems: (1) determining how to route an empty delivery vehicle back from its current location to its depot by a scheduled arrival time, and (2) selecting a profit-maximizing subset of delivery requests between various locations on the route subject to the vehicle's capacity. Figure \ref{fig:BPMP} illustrates a BPMP instance and solution.

Figure \ref{fig:BPMP} shows a network representation of the problem with an empty vehicle at a location represented by node 1. The vehicle weighs 1 ton and has a carrying capacity of $Q = 2$ tons of cargo. The vehicle needs to return to its depot, represented by node 6, within a fixed period of time. The vehicle’s average  traveling speed limits the route to node 6 to a maximum distance of 7 miles. The vehicle can make extra money by accepting delivery requests to pick up cargo at the locations represented by nodes 1 through 5, destined for locations represented by nodes 2 through 6 as long as it can get back to the depot on time. The tuple $(d_{ij}, w_{ij})$ indicates the distance (in miles) and the size of the delivery request (in tons) from node $i$ to node $j$. The optimal solution indicated in Figure 1 routes the vehicle on the path represented by the arc sequence (1, 3), (3, 5), (5, 6).

 BPMP was first introduced by Dong et al. (2006). Yu and Dong \cite{Yu} proposed a MIP formulation based on the traditional node-arc model of multicommodity flow. Dong \cite{Dong} proposed an alternative MIP formulation of BPMP called the triples formulation. Thus, in the literature there are two kinds of BPMP MIP formulations: node-arc and triples. The purpose of our experimental study is to enhance both models as much as possible by applying candidate techniques, and compare their performance using the CIM method.
 
 \subsection{Node-arc formulation}
 The following node-arc formulations is taken from \cite{BPMP}.  The binary variable $x_{ij}$ indicates whether or not the vehicle traverses arc $(i, j)$, and binary variable $y_{kl}$ indicates whether or not to accept request $(k, l)$. Binary variable $z_{kl,ij}$ determines whether or not request $(k, l)$ is performed via arc $(i, j)$. Variable $\theta_{ij}$ represents the total flow (i.e., tons of cargo) transported on arc $(i, j)$. Sequence variables $s_i \ge 0$, for $i= 1, \ldots, n$, track the relative order in which nodes are visited. The node-arc formulation for BPMP is
\begin{eqnarray}
\max_{s,x,y,z,\theta} & p \sum_{(k,l) \in \mathcal{R}} d_{kl} w_{kl} y_{kl} - c \sum_{(i, j) \in \mathcal{A}}  d_{ij} \theta_{ij} - c v \sum_{(i, j) \in \mathcal{A}} d_{ij} x_{ij}
\label{eqn:1}
\end{eqnarray}

\noindent subject to \begin{eqnarray}
  \sum_{j=2}^{n} x_{1j} & = & 1 \label{eqn:2e} \\
  \sum_{i = 1}^{n-1} x_{in} & = & 1 \label{eqn:2f} \\
  \sum_{i \in \mathcal{N} \setminus \{k,n\}} x_{ik} & = & \sum_{j \in \mathcal{N} \setminus \{1, k\}} x_{kj} \qquad \forall k \in \mathcal{N} \setminus \{1, n\} \label{eqn:2g}\\
  \sum_{(i,j) \in \mathcal{A}} d_{ij} x_{ij} & \le & D \label{eqn:2h} \\
  \sum_{i \in \mathcal{N} \setminus \{k, n\}} x_{ik} & \le & 1 \qquad \forall k \in \mathcal{N} \setminus \{1\} \label{eqn:node-degree} \\
  s_i - s_j + (n + 1) x_{ij} & \le & n   \qquad \forall (i, j) \in \mathcal{A} \label{eqn:2k} \\ \nonumber \\
  \sum_{(k,l) \in \mathcal{R}} z_{kl,ij} & \le & M x_{ij} \qquad \forall (i, j) \in \mathcal{A} \label{eqn:2d} \\
  \sum_{j \in \mathcal{N} \setminus \{1, k\}} z_{kl,kj} &  = &  y_{kl}   \qquad \forall (k, l) \in \mathcal{R} \label{eqn:2a} \\
  \sum_{i \in \mathcal{N} \setminus \{l,n\}} z_{kl,il} &  = &  y_{kl} \qquad \forall (k, l) \in \mathcal{R} \label{eqn:2b} \\
 \sum_{\{i \in \mathcal{N}: (i, h) \in \mathcal{A}\}} z_{kl,ih} & = & \sum_{\{j \in \mathcal{N}: (h, j) \in \mathcal{A}\}} z_{kl,hj}\qquad\forall (k, l) \in \mathcal{R}, h \in \mathcal{N} \setminus \{k,l\}           \label{eqn:2c} \\
  \theta_{ij} & = & \sum_{(k,l) \in \mathcal{R}} w_{kl} z_{kl,ij} \qquad (i, j) \in \mathcal{A} \label{eqn:2i} \\ \nonumber \\
  \theta_{ij} & \le & Q \qquad \qquad \forall (i, j) \in \mathcal{A} \label{eqn:2j} 
  \end{eqnarray}

The objective function is to maximize the total profit.The vehicle's route is constrained to at most $D$ miles by constraint (\ref{eqn:2h}). The {\em node-degree}  (\ref{eqn:node-degree}), and {\em subtour elimination} constraints (\ref{eqn:2k}) ensure that the vehicle follows a simple path from node 1 to node $n$. The sequence variables determine the relative order in which nodes are visited by the vehicle. The logical connection between $x_{ij}$ and $z_{kl,ij}$ is enforced by constraint set (\ref{eqn:2d}). Constraint sets (\ref{eqn:2a}) and (\ref{eqn:2b}) enforce the logical relationship between $y_{kl}$ and $z_{kl,ij}$.   Constraints  (\ref{eqn:2c}) are flow-conservation constraints for intermediate nodes on the path the vehicle takes from node $k$ to node $l$. The capacity limit is enforced by constraint set (\ref{eqn:2j}). 

 \subsection{Triples formulation}

\label{sec:TriplesMIP}
The following description of the triples formulation of BPMP is adapted from \cite{BPMP}; 
it uses a compact formulation of multicommodity flow originally proposed by Matula  \cite{NET:NET21583,Matula} in which  {\em triples variable} $u_{ij}^k$ for node triple $(i,j,k)$ represents the total flow on all paths from node $i$ to node $j$ with arc $(i,k)$ as the first arc. 
In the triples formulation of BPMP, $u_{ij}^k$ represents the tons of cargo that the vehicle carries from node $i$ to node $j$ on arc $(i,k)$ and an unspecified path from node $k$ to node $j$. In a feasible solution, the unspecified path turns out to be the route that the vehicle takes from node $k$ to node $j$ \cite{BPMP}.  The triples formulation of the BPMP replaces the $z$ variables of the node-arc formulation with triples variables. The multicommodity flow constraints  (\ref{eqn:2a})--(\ref{eqn:2i}) are replaced with the following set of {\em triples constraints} that relate the triples variables to arc flows:

\begin{equation}
\theta_{ij} = w_{ij} y_{ij} + \sum_{(i,k,j) \in \mathcal{T}} u_{ik}^j  + \sum_{(k,j,i) \in \mathcal{T}} u_{kj}^i -   \sum_{(i,j,k) \in \mathcal{T}} u_{ij}^k \qquad \forall (i, j) \in \mathcal{A}
\label{eqn:triples_constraints}
\end{equation}

\begin{equation}
  \theta_{ij}  \le  Q x_{ij} \qquad \forall (i,j) \in {\cal A}. \label{eqn:conditional arc flow}
\end{equation}
The following constraints are imposed in order to force arc $(i,k)$ to be on the vehicle's route
 if variable $u_{ij}^k$ is positive:

\begin{equation}
u_{ij}^k  \le  Q x_{ik} \qquad \forall (i, j, k) \in \mathcal{T}  
\label{eqn:link_u_x}
\end{equation}

\noindent These constraints provide a logical linkage between the $u$ variables and the $x$ variables, and replace constraint set (\ref{eqn:2d}) of the node-arc formulation. For a detailed explanation of the triples formulation, interested readers are referred to  \cite{BPMP}.

\section{Node-Arc Summary}
\label{sec:Node-Arc}
In this section we use the GCI to evaluate the efficacy of various techniques (cuts, branching rules, etc.) designed to improve CPLEX’s performance using the node-arc model given in Section \ref{sec:BPMP}.  These techniques were selected and informally ranked by effectiveness from a larger set of candidates after preliminary experiments that we performed prior to developing the CIM. Before applying the first of the techniques, we established an “incumbent” enhanced node-arc formulation by determining a tight Big-$M$ value for the $x$-$z$ linking constraint set (\ref{eqn:2d}).  We then applied the techniques sequentially according to the ranking from our preliminary experiments.
Using ``model 1” to refer to an incumbent solution approach and ``model 2” to refer to the application of a particular technique to model 1. 
If $\mbox{GCI}> 1$, we recommend adopting the technique, and making the resulting model the new incumbent. We also say, for convenience, that the model 2 is ``GCI times faster” than model 1. If $\mbox{GCI} \le 1$, we recommend not adopting the technique.

\subsection{Computing Environment and Weight Selection}

The computations were performed on Dell R730 computers each with Dual 12 Core 2.6GHz Intel Xeon processors and 380GB RAM.
The formulations were implemented in AMPL 10.00 and solved with CPLEX 12.6.0.0. We used the default settings for AMPL and CPLEX except where specified. 

In our experience, practitioners solving real-world problems are much more concerned about real time as a performance measure than CPU time, and are often unaware of the tick measure. For our purposes, however, the reproducibility of the tick metric is quite important. Therefore, we used the following weights for each type of time speedup: $\omega_C = 6, \omega_R = 8$, and $\omega_T = 8$. Thus, real time and ticks were equally important and more important than CPU time by a factor of $1 \frac{1}{3}$. For the ten instances of the same size problem, in reference to the PERT technique, we gave the largest weight to the median, the second largest to the mean, and the least weight to the min and max. For the problem sizes, we gave 30-node problems the largest weight, 20-node problems the second largest, and 10-node problems the least. The specific weights are listed in Table \ref{tab w}.
\begin{table}[h!]
\begin{center}
\caption{Node-Arc Weights}\label{tab w}
\begin{tabular}{|l|l|l|} \hline 
CPU: $\omega_c$ = 6         & Ticks: $\omega_t$ = 8         & Real Time: $\omega_t$ = 8  \\ \hline \hline
Median: $\omega_{med}$ = 40 & Mean:  $\omega_{mea}$ = 10    & Min \& Max: $\omega_{min}$ = $\omega_{max}$ = 0.5 \\ \hline \hline
10-node: $\omega_{10}$ = 1  & 20-node: $\omega_{20}$ = 10   & 30-node: $\omega_{30}$ = 12     \\ \hline
\end{tabular}
\end{center}
\end{table}

\subsection{Initial Incumbent Formulation}
We were only able to solve 10-node and 20-node instances with the incumbent solution approach (i.e., the original node-arc model).
The model was solved three times for each problem instance. 
With no techniques applied, the median average real time for the 10-node instances was 1.15 seconds, and the median average real time for the 20-node instances was about 20 minutes. The results for $n = 20$ are shown in Table \ref{tab2}. We summarize the results for $n = 10$ using four-tuples listing the averages of the min, mean, median, and max values for CPU time, real time, and ticks. The CPU seconds, real time seconds, and ticks tuples are (2.28, 4.17, 4.23, 6.45), (0.50, 1.15, 0.99, 2.56), and (120.06, 332.90, 272.06, 645.50), respectively. That is, the  maximum average CPU time for the 10-node instances was 6.45 seconds, the maximum average real time was 2.56 seconds, and the maximum average number of ticks was 645.50.
Hereinafter, CPU and real time are reported in seconds. 


\begin{table}
\begin{center}
\caption{Test Results of Original Node-Arc Model for $n$ = 20.}\label{tab2}
\begin{tabular}{|c|r|r|r|r|r|r|r|r|r|r|} \cline{2-9}
\multicolumn{1}{c}{}    & \multicolumn{4}{|c|}{CPU Time} & \multicolumn{4}{c|}{Real Time} & \multicolumn{1}{c}{} \\ \hline
Instance	&	Run 1	&	Run 2	&	 Run 3	&	Ave.	&	  Run 1	&	 Run 2	&	  Run 3	&	Ave. 	&	Ticks	\\ \hline \hline
1	        &	12,654	&	12,310	&	12,172	&	12,379	&	1,937	&	1,896	&	1,734	&	1,856	&		1,244,880	\\ \hline	
2	        &	1,745	&	1,757	&	1,716	&	1,739	&	355	    &  	372	    &	344	    &	357	    &		272,703	\\ \hline	
3	        &	2,457	&	2,544	&	2,393	&	2,465	&	497	    &	536	    &	467	    &	500	    &		382,350	\\ \hline	
4	        &	24,954	&	26,053	&	27,313	&	26,106	&	1,891	&	2,044	&	1,943	&	1,960	&		1,211,197	\\ \hline	
5	        &	2,964	&	3,036	&	2,981	&	2,994	&	632	    &	681	    &	628	    &	647	    &		583,158	\\ \hline	
6	        &	32,635	&	35,586	&	34,566	&	34,262	&	1,720	&	1,859	&	1,778	&	1,786	&		1,005,079	\\ \hline	
7	        &	2,760	&	2,770	&	2,615	&	2,715	&	502	    &	512	    &	482	    &	499	    &		422,516	\\ \hline	
8	        &	5,873	&	5,888	&	5,738	&	5,833	&	845	    &	857	    &	806	    &	836	    &		674,853	\\ \hline	
9	        &	28,620	&	32,115	&	28,746	&	29,827	&	2,007	&	2,226	&	1,967	&	2,067	&   	1,318,085	\\ \hline	
10	        &	17,096	&	17,745	&	18,202	&	17,681	&	1,571	&	1,749	&	1,606	&	1,642	&		1,091,225	\\ \hline \hline
Min	        &	1,745	&	1,757	&	1,716	&	1,739	&	355	    &	372	    &	344	    &	357	    &		272,703	\\ \hline	
Mean	    &	13,176	&	13,980	&	13,644	&	13,600	&	1,196	&	1,273	&	1,175	&	1,215	&		820,605	\\ \hline	
Median	    &	9,263	&	9,099	&	8,955	&	9,106	&	1,208	&	1,303	&	1,206	&	1,239	&		839,966	\\ \hline	
Max	        &	32,635	&	35,586	&	34,566	&	34,262	&	2,007	&	2,226	&	1,967	&	2,067	&		1,318,085	\\ \hline	
\end{tabular}
\end{center}
\end{table}

\subsection{Technique 1: Conditional Arc Flow}

The original node-arc model \cite{Yu} uses constraint set (\ref{eqn:2j}), $\theta_{ij} \le Q$, to ensure that the total amount of flow, $\theta_{ij}$, on arc $(i, j)$ is less than or equal to the vehicle capacity, $Q$. Notice that if the vehicle does not travel on arc $(i, j)$, there should be no flow on the arc (i.e., if $x_{ij}  = 0$, then $\theta_{ij} = 0$). If the vehicle does travel on arc $(i, j)$, the maximum flow on the arc is $Q$, (i.e., if $x_{ij} = 1$, then $\theta_{ij} \le Q$). Therefore, (\ref{eqn:2j}) can be replaced by the following constraint set which we call conditional arc-flow
\begin{equation}
\theta_{ij} \le Q x_{ij} \qquad \forall (i, j) \in \mathcal{A}. \label{eqn:conditional_arc_flow}
\end{equation}

Yu and Dong \cite{Yu} were unable to solve 30-node instances with the original node-arc model. We had a similar experience in our preliminary tests. Therefore, we tested this technique only on 10-, 20-, and 30-node instances. Table \ref{tab:cond_arc_flow_node_arc_20} gives detailed test results of three runs for the 20-node instances after applying the technique. The CPU time, real time, and ticks tuples for the 10-node instances are (0.90, 2.85, 2.91, 4.82), (0.26, 0.66, 0.61, 0.96), and (74.61, 237.37, 241.69, 384.78), respectively. The complete speedup summary is given in Table \ref{tab:cond_arc_flow_node_arc_incremental}. Table \ref{tab:cond_arc_flow_node_arc_GCI} lists the composite indices and GCI. Speedups in bold are greater than 1.

{\bf Conclusion}: After applying technique 1, conditional arc-flow, the GCI of speedups was 8.24, which means, on average, the model with conditional arc-flow was solved 8.24 times faster than the original model. Therefore, we adopted technique 1, replacing constraint set (\ref{eqn:2j}) with the conditional arc-flow constraints (\ref{eqn:conditional_arc_flow}). Furthermore, after applying conditional arc-flow constraints we were able to solve 30-node instances. The average real times ranged from 4,421 seconds (1.23 hours) to 17,829 seconds (4.95 hours) with a mean and median of 8,546 seconds (2.37 hours) and 7,635 seconds (2.12 hours), respectively \cite{BO}.

\begin{table}[ht!]
\begin{center}
\caption{Test Results of Incremental Effect of Conditional Arc-Flow for $n = 20$}\label{tab:cond_arc_flow_node_arc_20}
\begin{tabular}{|c|r|r|r|r|r|r|r|r|r|r|} \cline{2-9}
\multicolumn{1}{c}{}    & \multicolumn{4}{|c|}{CPU Time} & \multicolumn{4}{c|}{Real Time} & \multicolumn{1}{c}{} \\ \hline
Instance	    &	Run 1	    &	Run 2	    &	 Run 3	    &  	Ave.	&	  Run 1	&	 Run 2	&	  Run 3	&	Ave. 	&		Ticks	\\ \hline \hline
1	            &	1,306	&	1,258	&	1,266	&	1,276	&	96	&	100	&	97	&	98	    &		54,992	\\ \hline
2	            &	720	    &	701	    &	702	    &  	708	    &	63	&	63	&	64	&	63	    &		37,316	\\ \hline
3	            &	1,974	&	1,689	&	1,729	&	1,797	&	143	&	146	&	134	&	141	    &		81,335	\\ \hline
4	            &	4,945	&	3,848	&	3,836	&	4,210	&	358	&	362	&	318	&	346	    &		202,725	\\ \hline
5	            &	450	    &  	409	    &	408	    &	422	    &	94	&	94	&	92	&	93	    &		70,847	\\ \hline
6	            &	2,152	&	1,752	&	1,764	&	1,890	&	196	&	200	&	182	&	193	    &		119,363	\\ \hline
7	            &	1,243	&	1,105	&	1,120	&	1,156	&	81	&	81	&	77	&	80	    &		39,778	\\ \hline
8	            &	1,189	&	1,080	&	1,010	&	1,093	&	74	&	76	&	68	&	73	    &		35,113	\\ \hline
9	            &	673	    &	630	    &	620	    &	641	    &	127	&	130	&	125	&	128	    &	    96,582	\\ \hline
10	            &	811	    &	759	    &	740	    &	770	    &	123	&	125	&	119	&	122	    &		85,867	\\ \hline \hline
Min	            &	450	    &	409	    &	408	    &	422	    &	63	&	63	&	64	&	63	    &		35,113	\\ \hline
Mean	        &	1,546	&	1,323	&	1,320	&	1,396	&	136	&	138	&	128	&	134	    &		82,392	\\ \hline
Median	        &	1,216	&	1,093	&	1,065	&	1,125	&	110	&	112	&	108	&	110	    &		76,091	\\ \hline
Max	            &	4,945	&	3,848	&	3,836	&	4,210	&	358	&	362	&	318	&	346	    &		202,725	\\ \hline
\end{tabular}
\end{center}
\end{table}

\begin{table}[h!]
\begin{center}
\caption{Summary of Incremental Effect of Conditional Arc Flow Constraints. }\label{tab:cond_arc_flow_node_arc_incremental}
\begin{tabular}{|c|r|r|r|} \cline{2-4}
\multicolumn{1}{c}{} & \multicolumn{3}{|c|}{Speedup} \\ \hline
    	             &	Ave. CPU Time	&	Ticks 	&	Ave. Real Time	\\ \hline \hline
    	\multicolumn{4}{|c|}{$n = 10$} \\ \hline
Min 	            &	0.95	    &	0.91	    &	{\bf 1.10}	 \\ \hline
Mean	            &	{\bf 1.68}	&	{\bf 1.42}	&	{\bf 1.74}	 \\ \hline
Median	            &	{\bf 1.46}	&	{\bf 1.32}	&	{\bf 1.71}	 \\ \hline
Max	                &	{\bf 3.17}	&	{\bf 2.61}	&	{\bf 2.67}	\\ \hline 
\multicolumn{4}{|c|}{$n = 20$} \\ \hline
Min	                &	{\bf 1.37}	&	{\bf 4.70}	&   {\bf 3.55}  \\ \hline
Mean	            &	{\bf 12.21}	&	{\bf 11.35}	&	{\bf 9.74}	\\ \hline
Median	            &	{\bf 6.65}	&	{\bf 9.52}	&	{\bf 8.10}	\\ \hline
Max	                &	{\bf 46.52}	&	{\bf 22.64}	&	{\bf 19.01}	\\ \hline
\end{tabular}
\end{center}
\end{table}

\begin{table}[h!]
\begin{center}
\caption{Composite Indices and GCI for Technique 1}\label{tab:cond_arc_flow_node_arc_GCI}
\begin{tabular}{|c|c|c|c|c|c|} \hline
$n$	& $C_n$	(CPU)	        &	$T_n$ (Ticks)	    &	$R_n$ (Real Time)	&	$I_n$	    &	$\mbox{GCI}$	\\  \hline \hline
10	&	{\bf 1.52}	        &	{\bf 1.35}	        &	{\bf 1.72}	        &	{\bf 1.53}	&	{\bf 8.24}	\\ \cline{1-5}
20	&	{\bf 8.08}	        &	{\bf 9.96}	        &	{\bf 8.48}	        &	{\bf 8.91}	&		\\ \hline
\end{tabular}
\end{center}
\end{table}

\subsection{Technique 2: Relax Node-Degree Constraints}
Yu and Dong \cite{Yu} used the node-degree cuts (\ref{eqn:node-degree}) to ensure that the vehicle visits each location at most once. But, the MTZ subtour elimination constraints (\ref{eqn:2k}) also ensure that vehicle visits each node at most once in an integer solution. Therefore, we can relax (drop) the node-degree constraints without losing validity of the integer model (the node-degree cuts can be violated in solutions to the LP relaxations). Table \ref{tab8} gives the composite indices and GCI for this technique. Speedups greater than 1 are in bold. 

\begin{table}[h!]
\begin{center}
\caption{Composite Indices and GCI for Technique 2}\label{tab8}
\begin{tabular}{|c|c|c|c|c|c|} \hline
$n$	&	$C_n$ (CPU)	        &	$T_n$ (Ticks)	    &	$R_n$ (Real Time)	&	$I_n$	        &	$\mbox{GCI}$	    \\  \hline \hline
10	&	0.73	            &	{\bf 1.17}	        &	0.90	            &	{\bf 0.95}	    &		        \\ \cline{1-5}
20	&	{\bf 1.02}	        &	{\bf 1.25}	        &	{\bf 1.22}	        &	{\bf 1.17}      &	{\bf 1.28}	\\ \cline{1-5}
30	&	{\bf 1.43}          &	{\bf 1.28}	        &	{\bf 1.48}	        &	{\bf 1.40}	    &               \\ \hline
\end{tabular}
\end{center}
\end{table}

{\bf Conclusion:} After applying technique 2, relax node-degree constraints, the CGI of speedups was 1.28, which means, on average, the model relaxing the node-degree constraints was solved 1.28 times faster than the incumbent model. Therefore, we adopted technique 2 and dropped constraint set (\ref{eqn:node-degree}) from the incumbent. 

\subsection{Technique 3: Single-Node Demand Cuts}
The single-node demand cuts state that the total weight of the delivery requests accepted from node $i$, or into node $j$, is at most the vehicle capacity, $Q$. 

\begin{eqnarray}
&  \sum_{j \in V \setminus \left\{1,i\right\}}{w_{ij}y_{ij}}\le Q &  \forall i\in V\setminus\left\{n\right\} \\
&  \sum_{i\in V \setminus \left\{j,n\right\}} {w_{ij}y_{ij}}\le Q &  \forall j\in V\setminus\left\{1\right\}
\end{eqnarray}

The above are valid inequalities that are satisfied by any feasible solution because the vehicle cannot simultaneously hold cargoes with total weight more than its capacity. This condition is not necessarily enforced by solutions to the LP relaxation because of the fractional $y$ values. Table  \ref{tab9} gives the results from applying single-node demand cuts to the incumbent node-arc model.

\begin{table}[h!]
\begin{center}
\caption{Composite Indices and GCI Technique 3}\label{tab9}
\begin{tabular}{|c|c|c|c|c|c|} \hline
$n$	&	$C_n$ (CPU)	        &	$T_n$ (Ticks)	    &	$R_n$ (Real Time)	&	$I_n$	        &	$\mbox{GCI}$	    \\  \hline \hline
10	&	{\bf 1.30} 	&	0.96	&	{\bf 1.21}	&	{\bf 1.15}	&		        \\ \cline{1-5}
20	&	    0.86	&	0.91	&	0.88    	&	0.89        &	0.94	\\ \cline{1-5}
30	&	{\bf 1.01}  &	0.95	&	0.94	    &	0.96	    &               \\ \hline
\end{tabular}
\end{center}
\end{table}

{\bf Conclusion:} After applying technique 3, single-node demand cuts, the GCI of speedups was 0.94, which means that solving the incumbent model was faster. Therefore, we did not adopt technique 3. 

\subsection{Technique 4: Relax $x$-$z$ Linking Constraints}
Adopting the conditional arc-flow cuts makes the constraints linking the $x$ and $z$ variables  redundant.
Therefore, we can relax (\ref{eqn:2d}).  Table  \ref{tab10} gives the results from applying relax $x$-$z$ linking constraints to the incumbent node-arc model.

\begin{table}[h!]
\begin{center}
\caption{Composite Indices and GCI Technique 4}\label{tab10}
\begin{tabular}{|c|c|c|c|c|c|} \hline
$n$	&	$C_n$ (CPU)	        &	$T_n$ (Ticks)	    &	$R_n$ (Real Time)	&	$I_n$	        &	$\mbox{GCI}$	    \\  \hline \hline
10	&	{\bf 1.50} 	&	{\bf 1.13}	&	{\bf 1.16}	&	{\bf 1.24}	&		        \\ \cline{1-5}
20	&	    {\bf 1.06}	&	{\bf 1.12}	&	{\bf 1.08}    	&	{\bf 1.09}        &	{\bf 1.18}	\\ \cline{1-5}
30	&	{\bf 1.39}  &	{\bf 1.16}	&	{\bf 1.21}	    &	{\bf 1.24}	    &               \\ \hline
\end{tabular}
\end{center}
\end{table}

{\bf Conclusion:} After applying technique 4, relax $x$-$z$ linking constraints, the grand composite index of speedups (GCI) was 1.18, which means, on average, solving the model relaxing $x$-$z$ linking constraints was 1.18 times faster than the incumbent model. Therefore, we adopted technique 4 and dropped constraint set (\ref{eqn:2d}) from the incumbent. 

\subsection{Technique 5: Branching Priority}
In the node-arc model, there are three types of binary variables: $x_{ij}$, $y_{kl}$ and $z_{kl,ij}$ which indicate whether the vehicle travels on arc $(i, j)$, whether the delivery request from node $k$ to node $l$ is accepted, and whether the accepted demand from node $k$ to node $l$ is realized via arc $(i, j)$, respectively. We suspected that prioritizing determining the vehicle’s route over deciding which delivery requests to accept would lead to faster solution times. Therefore, we tested solving the problem with a branching rule stating that $x$ variables are branched on before any other binary variables.  
 Table  \ref{tab11} gives the results from applying branching priority to the incumbent node-arc model.

\begin{table}[h!]
\begin{center}
\caption{Composite Indices and GCI Technique 5}\label{tab11}
\begin{tabular}{|c|c|c|c|c|c|} \hline
$n$	&	$C_n$ (CPU)	        &	$T_n$ (Ticks)	    &	$R_n$ (Real Time)	&	$I_n$	        &	$\mbox{GCI}$	    \\  \hline \hline
10	&	0.75 	&	1.0	&	{\bf 1.05}	&	0.95	&		        \\ \cline{1-5}
20	&	    {\bf 1.86}	&	{\bf 1.50}	&	{\bf 1.42}    	&	{\bf 1.57}        &	{\bf 1.53}	\\ \cline{1-5}
30	&	{\bf 1.81}  &	{\bf 1.36}	&	{\bf 1.51}	    &	{\bf 1.54}	    &               \\ \hline
\end{tabular}
\end{center}
\end{table}

{\bf Conclusion:} After applying technique 5, branching priority, the grand composite index of speedups (GCI) was 1.5264, which means that using the branching rule was an improvement over solving the incumbent model with CPLEX’s default settings. Therefore, we adopted the branching priority on the $x$ variables. Hereinafter we refer to the process of solving the incumbent model with the branching rule as the “incumbent model”.

\subsection{Technique 6: Lifted MTZ }

Desrochers and Laporte \cite{desrochers1991improvements} proved that the MTZ subtour elimination constraints (\ref{eqn:2k}) can be strengthened by lifting them to
\begin{equation}
      s_i - s_j + (n - 1) x_{ij} + (n-3) x_{ji} \le  n - 2  \qquad \forall  i \in {\mathcal N} \setminus \{1,n\}, j \in {\mathcal N \setminus \{1, i, n\}}. \label{eqn:lifted mtz}
\end{equation}

 \noindent Table  \ref{tab12} gives the results from applying lifted MTZ to the incumbent node-arc model.

\begin{table}[h!]
\begin{center}
\caption{Composite Indices and GCI Technique 6}\label{tab12}
\begin{tabular}{|c|c|c|c|c|c|} \hline
$n$	&	$C_n$ (CPU)	        &	$T_n$ (Ticks)	    &	$R_n$ (Real Time)	&	$I_n$	        &	$\mbox{GCI}$	    \\  \hline \hline
10	&	{\bf 1.88} 	&	{\bf 1.04} 	&	{\bf 1.14}	&	{\bf 1.30}	&		        \\ \cline{1-5}
20	&	    0.98	&	0.99	&	0.99    	&	0.99        &	0.9991	\\ \cline{1-5}
30	&	{\bf 1.22}  &	0.82	&	0.97	    &	0.98	    &               \\ \hline
\end{tabular}
\end{center}
\end{table}

{\bf Conclusion:} After applying technique 6, lifted MTZ, the grand composite index of speedups (GCI) was 0.9991, which means that the incumbent model was solved faster. Therefore, we did not adopt technique 6. 

\subsection{Technique 7: MTZ upper bound}
In the original MTZ subtour elimination constraints \cite{Miller}, there is no upper limit for the sequence variable $s_i$. As a result any given tour has essentially an infinite number of representations in terms of the sequence variables. This type of symmetry can needlessly slow down the branch-and-bound process by causing it “to explore and eliminate such alternative symmetric solutions” \cite{sherali2001improving}. 
Desrochers and Laporte \cite{desrochers1991improvements} proved that following constraints ensure a unique representation of any given feasible tour:

\begin{equation}
1 \le s_i \le (n-1) \qquad \forall i \in {\mathcal N} \setminus \{1\}. \label{eqn:mtz upper bound}
\end{equation}

\noindent Table  \ref{tab13} gives the results from applying MTZ upper bound to the incumbent node-arc model.

\begin{table}[h!]
\begin{center}
\caption{Composite Indices and GCI Technique 7}\label{tab13}
\begin{tabular}{|c|c|c|c|c|c|} \hline
$n$	&	$C_n$ (CPU)	        &	$T_n$ (Ticks)	    &	$R_n$ (Real Time)	&	$I_n$	        &	$\mbox{GCI}$	    \\  \hline \hline
10	&	{\bf 1.49} 	&	0.97	&	0.99	&	{\bf 1.12}	&		        \\ \cline{1-5}
20	&	    0.98	&	0.98	&	0.94    	&	0.97        &	0.9096	\\ \cline{1-5}
30	&	0.85  &	0.93	&	0.76	    &	0.84	    &               \\ \hline
\end{tabular}
\end{center}
\end{table}

{\bf Conclusion:} After applying technique 7, MTZ upper bound, the grand composite index of speedups (GCI) was 0.9096 indicating that it was more efficient to solve the incumbent model.  Therefore, we did not adopt the upper bound constraints for the MTZ sequence variables.  
\subsection{Technique 8: Cover Cuts}
Fischetti et al. \cite{fischetti1998solving} found that cover cuts on sets of arcs whose total length is more than the maximum route distance $D$ were effective for solving the Orienteering Problem, which is a special case of the BPMP.  They also proposed solving a knapsack problem to determine if there is a set of arcs $\mathcal{S}$ that violates the cover cut:

\begin{equation}
\sum_{(i, j) \in {\mathcal S} } x_{ij} \le |\mathcal{S}|-1 \qquad  \forall \mathcal{S} \subseteq \mathcal{A} \mbox{ such that } \sum_{(i,j) \in \mathcal{S}} d_{ij} > D \label{cover cut}
\end{equation}

\noindent in the LP relaxation.  We applied this technique iteratively to the BPMP adding violated cover cuts as necessary until no additional cover cuts are found at which point we solved the MIP. Table  \ref{tab14} gives the results from applying cover cuts to the incumbent node-arc model.

\begin{table}[h!]
\begin{center}
\caption{Composite Indices and GCI Technique 8}\label{tab14}
\begin{tabular}{|c|c|c|c|c|c|} \hline
$n$	&	$C_n$ (CPU)	        &	$T_n$ (Ticks)	    &	$R_n$ (Real Time)	&	$I_n$	        &	$\mbox{GCI}$	    \\  \hline \hline
10	&	{\bf 1.20} 	&	0.97	&	0.74	&	0.95	&		        \\ \cline{1-5}
20	&	    0.79	&	0.98	&	0.77    	&	0.85        &	0.9069	\\ \cline{1-5}
30	&	{\bf 1.04}  &	0.96	&	0.88	    &	0.95	    &               \\ \hline
\end{tabular}
\end{center}
\end{table}

{\bf Conclusion:} After applying technique 8, cover cuts, the grand composite index of speedups (GCI) is 0.9069 indicating that it was more efficient to solve the incumbent model.  Therefore, we did not adopt the technique of cover cuts.  
\subsection{Technique 9: Pairwise Demand Cuts}

Pairwise demand cuts state that pairs of delivery requests from the same node whose total weight exceeds the vehicle's capacity are mutually exclusive:
\begin{equation}
      y_{ki} + y_{kj} \le  1 \quad \forall \{(k, i), (k, j) \in  A: w_{kl}  + w_{kj} > Q , i \neq j\} .  \label{eqn:pairwise demand cuts}
\end{equation}

\noindent The above are valid inequalities that are satisfied by any feasible solution to the MIP formulation. 
There can be a relatively large number of these cuts. Instead of adding all of them to the node-arc model, we adopt a simple scheme to add them as necessary. That is, we check for violated pairwise demand cuts of the corresponding LP relaxation, add any violated cuts found to the model, and solve the LP again. This process is repeated until no more cuts are found in the LP relaxation problem, at which point we restore the integrality constraints and solve the MIP. In this way, we can use a minimal number of pairwise demand cuts.    Table  \ref{tab15} gives the results from applying pairwise demand cuts to the incumbent node-arc model.

\begin{table}[h!]
\begin{center}
\caption{Composite Indices and GCI Technique 9}\label{tab15}
\begin{tabular}{|c|c|c|c|c|c|} \hline
$n$	&	$C_n$ (CPU)	        &	$T_n$ (Ticks)	    &	$R_n$ (Real Time)	&	$I_n$	        &	$\mbox{GCI}$	    \\  \hline \hline
10	&	{\bf 1.65} 	&	0.97	&	0.74	&	{\bf 1.07}	&		        \\ \cline{1-5}
20	&	    0.83	&	0.97	&	0.80    	&	0.87        &	0.9853	\\ \cline{1-5}
30	&	{\bf 1.11}  &	{\bf 1.05}	&	{\bf 1.07}	    &	{\bf 1.08}	    &               \\ \hline
\end{tabular}
\end{center}
\end{table}

{\bf Conclusion:} After applying technique 9, pairwise demand cuts, the grand composite index of speedups (GCI) was 0.9853, which means that the incumbent model was solved faster. Therefore, we did not adopt technique 9. 

\subsection{Best Node-Arc Model}
In total, we tested a series of nine techniques for improving solution time using CPLEX to solve the node-arc formulation of BPMP. Four of the techniques were adopted resulting in the Best Node-Arc Model \cite{BO}. Table \ref{best node arc vs original node arc} summarizes the speedup of the Best Node-Arc Model compared to the original model proposed by Yu and Dong \cite{Yu}. 
The CPU time, real time, and ticks tuples for $n = 30$ are
(20,454, 55,459, 39,485, 140,354), (1,594, 3,198, 2,564, 6,166), and (812,789, 1,192,684, 1,147,441, 1,966,885), respectively. The CPU time, real time, and ticks tuples for $n = 40$ are (334,025, 6,643,337, 1,213,615, 52,502,367),  (18,413, 329,773, 56,428, 2,652,518), and (6,601,682, 62,036,444, 15,873,253, 463,811,772), respectively.

\begin{table}[h!]
\begin{center}
\caption{Best Node-Arc Model vs. Original Node-Arc Model} \label{best node arc vs original node arc}
\begin{tabular}{|c|r|r|r||r|r|r|} \cline{2-7}
\multicolumn{1}{c|}{} & \multicolumn{6}{c|}{Speedup} \\ \cline{2-7}
\multicolumn{1}{c|}{} & \multicolumn{3}{c||}{$n = 10$} & \multicolumn{3}{c|}{$n = 20$} \\ \cline{2-7}
\multicolumn{1}{c|}{} & CPU & Ticks & Real Time& CPU & Ticks & Real Time \\ \hline
Min	                  &	0.58	&	{\bf 1.26}	&	{\bf 1.29}	&	{\bf 2.17}	&	{\bf 12.38}	&	{\bf 7.76}	\\ \hline
Mean	              &	{\bf 1.24}	&	{\bf 1.88}	&	{\bf 1.81}	&	{\bf 26.24}	&	{\bf 25.61}	&	{\bf 20.78}	\\ \hline
Median	              &	0.89	&	{\bf 1.9}	&	{\bf 1.63}	&	{\bf 17.97}	&	{\bf 24.58}	&	{\bf 19.56}	\\ \hline
Max	                  &	{\bf 3.68}	&	{\bf 2.86}	&	{\bf 2.81}	&	{\bf 111.29}	&	{\bf 42.77}	&	{\bf 43.98}	\\ \hline
\end{tabular}
\end{center}
\end{table}


\section{Triples Summary}
\label{sec:Triples}
In \cite{BO} we described in detail how we applied the CIM to the triples formulation of the BPMP proposed by Dong \cite{Dong}. We tested two techniques that are specific to the triples formulation, and six of the nine techniques we tested for the node-arc model.  Due to the fact that preliminary studies showed that we could solve larger problem instances with the triples formulation than with the node-arc formulation, we used the following instance-size weights: $\omega_{10} = 6$,  $\omega_{20} = 10$, $\omega_{30} = 13$,  $\omega_{40} = 14$, and  $\omega_{50} = 16$.  The other weights were the same as those used for the node-arc formulation.

\subsection{Technique 10: Relax Triples Linking Constraints}

The linking constraint (\ref{eqn:link_u_x}) forces ($i$, $k$) to be an arc on the vehicle's route if variable $u_{ij}^k$  is positive. However, Dong \cite{Dong} showed that model remains valid even if this constraint is relaxed. Relaxing (\ref{eqn:link_u_x}) significantly reduces the number of constraints in the triples model and consequently improves solution time.  Table  \ref{tab16} gives the results from applying relax triples linking constraints to the incumbent triples model.

\begin{table}[h!]
\begin{center}
\caption{Composite Indices and GCI Technique 10}\label{tab16}
\begin{tabular}{|c|c|c|c|c|c|} \hline
$n$	&	$C_n$ (CPU)	        &	$T_n$ (Ticks)	    &	$R_n$ (Real Time)	&	$I_n$	        &	$\mbox{GCI}$	    \\  \hline \hline
10	&	0.89 	&	{\bf 1.17}	&	{\bf 1.37}	&	{\bf 1.17}	&		        \\ \cline{1-5}
20	&	    {\bf 1.37}	&	{\bf 1.41}	&	{\bf 1.53}    	&	{\bf 1.44} 	&		        \\ \cline{1-5}       
30	&	{\bf 1.11}  &	{\bf 1.05}	&	{\bf 1.07}	    &	{\bf 1.08}	   &	{\bf 2.11}	\\ \cline{1-5}
40	&	    {\bf 3.05}	&	{\bf 1.55}	&	{\bf 3.36}    	&	{\bf 2.62} 	&		        \\ \cline{1-5}
50	&	    {\bf 1.64}	&	{\bf 1.09}	&	{\bf 1.93}    	&	{\bf 1.55} &               \\ \hline
\end{tabular}
\end{center}
\end{table}

Conclusion: After applying technique 10, relax triples linking constraints, the grand composite index of speedups (GCI) was 2.11, which means, on average, the triples model with (16) relaxed was solved 2.11 times faster than the incumbent model. Thus, we adopted it. 

\subsection{Technique 11: Enforce Node-Degree }

Unlike the node-arc model, the original triples model does not explicitly enforce the node-degree constraints (\ref{eqn:node-degree})  because the MTZ subtour elimination constraints (\ref{eqn:2k}) ensure that vehicle visits each node at most once in an integer solution. However, the node-degree constraints can be violated in solutions to the LP relaxation of the triples model.   Table  \ref{tab17} gives the results from applying enforce node-degree to the incumbent triples model.

\begin{table}[h!]
\begin{center}
\caption{Composite Indices and GCI Technique 11}\label{tab17}
\begin{tabular}{|c|c|c|c|c|c|} \hline
$n$	&	$C_n$ (CPU)	        &	$T_n$ (Ticks)	    &	$R_n$ (Real Time)	&	$I_n$	        &	$\mbox{GCI}$	    \\  \hline \hline
10	&	{\bf 1.27} 	&	{\bf 1.04}	&	{\bf 1.18}	&	{\bf 1.15}	&		        \\ \cline{1-5}
20	&	    {\bf 1.14}	&	0.98	&	{\bf 1.06}    	&	{\bf 1.05} 	&		        \\ \cline{1-5}       
30	&	{\bf 1.54}  &	{\bf 1.31}	&	{\bf 1.29}	    &	{\bf 1.37}	   &	{\bf 1.43}	\\ \cline{1-5}
40	&	    {\bf 1.54}	&	{\bf 1.38}	&	{\bf 1.45}    	&	{\bf 1.45} 	&		        \\ \cline{1-5}
50	&	    {\bf 2.13}	&	{\bf 1.80}	&	{\bf 1.56}    	&	{\bf 1.80} &               \\ \hline
\end{tabular}
\end{center}
\end{table}

Conclusion: After applying technique 11, enforce node-degree constraints, the grand composite index of speedups (GCI) was 1.43, which means, on average, the model with enforce node-degree constraints was solved 1.43 times faster than the incumbent model. Thus, we adopted it. 

\subsection{Best Triples Model}
Tables \ref{tab:best triples vs original triples 1} and \ref{tab:best triples vs original triples 2} summarize the speedup of the Best Triples Model compared to the original Triples model. Using the Best Triples Model, we were able to solve 50-node instances with a maximum average real time of 5,016 seconds (1.4 hours) \cite{BO,BPMP}.  Tables \ref{tab:best triples vs best node arc 1} and \ref{tab:best triples vs best node arc 2} compare the final triples and node-arc formulations. We did not calculate the GGI for Tables \ref{tab:best triples vs best node arc 1} and  \ref{tab:best triples vs best node arc 2}, but it is clearly larger than 1. Thus, the Best Triples Model is the solution approach recommended by the CIM.

\begin{table}[h!]
\begin{center}
\caption{Best Triples vs. Original Triples: $n = 10$, $n = 20$, and $n = 30$}
\label{tab:best triples vs original triples 1}
\begin{tabular}{|c|r|r|r||r|r|r||r|r|r|} \cline{2-10}
\multicolumn{1}{c|}{} & \multicolumn{9}{c|}{Speedup} \\\cline{2-10}
\multicolumn{1}{c|}{} & \multicolumn{3}{c||}{$n = 10$} & \multicolumn{3}{c||}{$n = 20$} & \multicolumn{3}{c|}{$n = 30$} \\ \cline{2-10}
\multicolumn{1}{c|}{} & CPU & Ticks & Real & CPU & Ticks & Real & CPU & Ticks & Real \\ \hline
Min	    &	{\bf 1.01}	&	0.49	    &	0.44	    &	0.80	    &	0.74	    &	0.84	    &	0.71	    &	0.69	    &	0.60	\\ \hline	
Mean	&	{\bf 1.51}	&	{\bf 1.11}	&	{\bf 1.88}	&	{\bf 1.95}	&	{\bf 1.66}	&	{\bf 1.69}	&	{\bf 6.49}	&	{\bf 3.13}	&	{\bf 4.48}	\\ \hline	
Median	&	{\bf 1.48}	&	{\bf 1.02}	&	{\bf 1.87}	&	{\bf 1.53}	&	{\bf 1.17}	&	{\bf 1.28}	&	{\bf 6.43}	&	{\bf 3.18}	&	{\bf 4.50}	\\ \hline	
Max	    &	{\bf 2.76}	&	{\bf 2.28}	&	{\bf 3.20}	&	{\bf 6.20}	&	{\bf 6.14}	&	{\bf 5.43}	&	{\bf 12.03}	&	{\bf 4.70}	&	{\bf 7.89}	\\ \hline	
\end{tabular}
\end{center}
\end{table}

\begin{table}[h!]
\begin{center}
\caption{Best Triples vs. Original Triples: $n = 40$ and $n = 50$}
\label{tab:best triples vs original triples 2}
\begin{tabular}{|c|r|r|r||r|r|r|} \cline{2-7}
\multicolumn{1}{c|}{} & \multicolumn{6}{c|}{Speedup} \\\cline{2-7}
\multicolumn{1}{c|}{} & \multicolumn{3}{c||}{$n = 40$} & \multicolumn{3}{c|}{$n = 50$} \\ \cline{2-7}
\multicolumn{1}{c|}{} & CPU & Ticks & Real & CPU & Ticks & Real\\ \hline
Min	    &	{\bf 1.49}	&	{\bf 1.11}	&	{\bf 1.31}	&	{\bf 3.91}	&	{\bf 2.00}	&	{\bf 4.62}	\\ \hline	
Mean	&	{\bf 7.97}	&	{\bf 3.17}	&	{\bf 6.59}	&	{\bf 10.34}	&	{\bf 3.65}	&	{\bf 9.07}	\\ \hline	
Median	&	{\bf 7.74}	&	{\bf 3.40}	&	{\bf 7.37}	&	{\bf 8.53}	&	{\bf 3.54}	&	{\bf 9.24}	\\ \hline	
Max	    &	{\bf 14.87}	&	{\bf 5.30}	&	{\bf 12.27}	&	{\bf 24.11}	&	{\bf 5.82}	&	{\bf 17.17}	\\ \hline	
\end{tabular}
\end{center}
\end{table}

\begin{table}[h!]
\begin{center}
\caption{Best Triples vs. Best Node Arc: $n = 10$ and $n = 20$}
\label{tab:best triples vs best node arc 1}
\begin{tabular}{|c|r|r|r||r|r|r|} \cline{2-7}
\multicolumn{1}{c|}{} & \multicolumn{6}{c|}{Speedup} \\\cline{2-7}
\multicolumn{1}{c|}{} & \multicolumn{3}{c||}{$n = 10$} & \multicolumn{3}{c|}{$n = 20$}  \\ \cline{2-7}
\multicolumn{1}{c|}{} & CPU & Ticks & Real & CPU & Ticks & Real  \\ \hline
Min	&	0.39	&	{\bf 1.47}	&	0.86	&	{\bf 6.34}	&	{\bf 8.89}	&	{\bf 7.98}	\\ \hline	
Mean	&	{\bf 2.38}	&	{\bf 4.47}	&	{\bf 2.31}	&	{\bf 62.43}	&	{\bf 24.02}	&	{\bf 21.46}	\\ \hline	
Median	&	{\bf 2.4}	&	{\bf 4.08}	&	{\bf 2.22}	&	{\bf 68.97}	&	{\bf 25.39}	&	{\bf 20.01}	\\ \hline	
Max	&	{\bf 4.50}	&	{\bf 8.38}	&	{\bf 4.00}	&	{\bf 137.06}	&	{\bf 42.75}	&	{\bf 46.20}	\\ \hline	
\end{tabular}
\end{center}
\end{table}

\begin{table}[h!]
\begin{center}
\caption{Best Triples vs. Best Node Arc: $n = 30$ and $n = 40$}
\label{tab:best triples vs best node arc 2}
\begin{tabular}{|c|r|r|r||r|r|r|} \cline{2-7}
\multicolumn{1}{c|}{} & \multicolumn{6}{c|}{Speedup} \\\cline{2-7}
\multicolumn{1}{c|}{} & \multicolumn{3}{c||}{$n = 30$} & \multicolumn{3}{c|}{$n = 40$}  \\ \cline{2-7}
\multicolumn{1}{c|}{} & CPU & Ticks             & Real              & CPU               & Ticks             & Real  \\ \hline
Min	    &	{\bf 34.33}	    &	{\bf 16.06}	    &	{\bf 30.96}	    &	{\bf 44.46}	    &	{\bf 34.79}	    &	{\bf 56.61}	\\ \hline
Mean	&	{\bf 210.55}	&	{\bf 52.26}	    &	{\bf 101.59}	&	{\bf 1,943.68}	&	{\bf 319.00}	&	{\bf 1,050.75}	\\ \hline
Median	&	{\bf 153.62}	&	{\bf 44.40}	    &	{\bf 106.01}	&	{\bf 327.41}	&	{\bf 96.47}	    &	{\bf 225.86}	\\ \hline
Max	    &	{\bf 872.08}	&	{\bf 137.89}	&	{\bf 212.53}	&	{\bf 11,929.94}	&	{\bf 1,417.70}	&	{\bf 4,275.04}	\\ \hline
\end{tabular}
\end{center}
\end{table}

\section{Conclusions}
\label{sec:Conclusion}
We have shown that the Composite Index Method (CIM) fills a gap in the area of optimization benchmarking. Calculating a single index, GCI, makes it much easier to select the best solution approach among multiple candidates.  By applying CIM to the backhaul profit maximization problem (BPMP), we demonstrate the step-by-step details of the framework of CIM in a parallel computing environment. Although we focused on solution-time measures for finding a provably optimal solution, the CIM can be easily adapted to consider other dimensions of concern such as memory usage and solution quality (for heuristics).
Furthermore, individual users can use their own weighting scheme to emphasize their personal preferences for making trade-offs between performance measures.

 Although the Best Triples Model of the BPMP that we identified by our iterative application of the CIM is a significant improvement over the initial node-arc formulation proposed in the literature, it is possible that we could have discovered an even better model by testing the model enhancements (techniques) in a different order. Thus, an important question for future research is how to determine the order in which alternative solution approaches are compared using the CIM when the approaches are not mutually exclusive.
 
 The successful utilization of the CIM lies in the wise selection of performance measures and weights in each dimension of concern. An illustrative study applying the CIM to well known optimization problems for different use cases is planned for the future.

%
%
%

\bibliographystyle{splncs04}
\bibliography{LION}

\end{document}